\documentclass[11pt,reqno]{amsproc}
\usepackage{natbib,xcolor}
\usepackage[all]{xy}
\usepackage[colorlinks,urlcolor=blue,citecolor=blue,linkcolor=blue]{hyperref}

\usepackage{mathpazo,graphicx,bm,anysize,caption,subcaption,eqnarray}

\marginsize{15.3mm}{15.3mm}{6.3mm}{5.25mm}

\theoremstyle{theorem}
\newtheorem{theorem}{Theorem}[section]

\theoremstyle{remark}

\theoremstyle{definition}

\newtheorem{example}[theorem]{Example}

\AtBeginDocument{\ \bigskip}

\makeatletter
\renewcommand{\bibsection}{\@bibtitlestyle}
\makeatother

\numberwithin{equation}{section}
\allowdisplaybreaks

\newcommand{\cdummy}{\cdot}
\newcommand{\mathd}{\mathrm{d}}

\begin{document}
\title{On the integral formula of the Jacobian determinant}

\author{Shibo Liu\vspace{-1em}}
\dedicatory{Department of Mathematics \& Systems Engineering, Florida Institute of Technology\\
Melbourne, FL 32901, USA}
\thanks{\emph{Emails}: \texttt{\bfseries sliu@fit.edu} (S. Liu)}

\begin{abstract}
  It is known that the integral of the Jacobian determinant of a smooth map $f
  : \bar{\Omega} \rightarrow \mathbb{R}^n$ depends only on $f |_{\partial
  \Omega}  $ and this result leads to an analytic proof of the
  Brouwer fixed point theorem. In this note we provide two new proofs of this
  result, one by classical analysis and one by differential forms and Stokes
  formula.
\end{abstract}

\maketitle

\section{Introduction}

In a recent paper {\cite{MR4725395}}, to respond inquiry from some readers
about the solution of Exercise 1.2.1 in his book {\cite{MR2435520}}, Krylov
gives a proof of the following theorem, see \citet[Lemma 1]{MR4725395}, where $f_\pm$ only need to be $C^1$.

\begin{theorem}
  \label{l1}Let $\Omega$ be a bounded domain in $\mathbb{R}^n$, $f_{\pm} :
  \bar{\Omega} \rightarrow \mathbb{R}^n$ be $C^2$-maps such that $f_+
  |_{\partial \Omega}   = f_- |_{\partial \Omega}  $. then
  \begin{equation}
    \int_{\Omega} \det f_+' (x) \mathd x = \int_{\Omega} \det f_-' (x) \mathd
    x \text{.} \label{e1}
  \end{equation}
\end{theorem}
%\end{document}
Using this theorem, the \emph{no-retraction theorem} follows immediately, see \citet[Corollary
1]{MR4725395}. Let $B$ be the closed unit ball in
$\mathbb{R}^n$. If there was a smooth retraction $R : B \rightarrow \partial
B$, since $R |_{\partial B}   = \mathrm{id} |_{\partial B}  $,
Theorem \ref{l1} with $f_+ = R$, $f_- = \mathrm{id}$ yields a contradiction
\begin{equation}
  \int_B \det R' (x) \mathd x = \int_B \mathd x = \mathrm{Vol} (B) > 0 \text{,}
  \label{e2}
\end{equation}
because $\det R' (x) = 0$, a consequence of $R (B) \subset \partial B$. That the two
integrals in (\ref{e2}) are equal also follows from a version of the changing
variable formula given in {\citet[Theorem 3.1]{MR3711061}}: \emph{For smooth closed
bounded domain $D$ in $\mathbb{R}^n$ and smooth map $\varphi : B \rightarrow
D$, if $\varphi : \partial B \rightarrow \partial D$ is a diffeomorphism, then
for continuous $f : D \rightarrow \mathbb{R}$ there holds}
\[ \int_D f (y) \mathd y = \pm \int_B f (\varphi (x)) \det \varphi' (x) \mathd
   x \text{.} \]
The first equality in (\ref{e2}) follows by letting $f = 1$, $D = B$ and
$\varphi = R$ in this formula. It is well known that the no-retraction theorem
is equivalent to the famous Brouwer fixed point theorem.

Krylov's proof of Theorem \ref{l1} in {\cite{MR4725395}} is based on the
observation: for small $t$, $f_{\pm}^t = \mathrm{id} + t f_{\pm}$ are
diffeomorphisms and $f^t_+ (\Omega) = f_-^t (\Omega)$ because $f_+ |_{\partial
\Omega}   = f_- |_{\partial \Omega}  $, thus by the changing
variable formula
\[ \int_{\Omega} \det (I + t f_+' (x)) \mathd x = \mathrm{Vol} (f_+^t (\Omega))
   = \mathrm{Vol} (f_-^t (\Omega)) = \int_{\Omega} \det (I + t f_-' (x)) \mathd
   x \text{.} \]
Since the two sides are polynomials in $t$, comparing the coefficients of
$t^n$ gives (\ref{e1}). To make the argument rigorous some issues including
why $f_{\pm}^t (\partial \Omega) = \partial f_{\pm}^t (\Omega)$ and why
$\partial f_+^t (\Omega) = \partial f_-^t (\Omega)$ implies $f_+^t (\Omega) =
f_-^t (\Omega)$ need to be handled. These are clarified via a careful point
set analysis in a long remark following the proof. One of the advantages of
this argument is that the maps $f_{\pm}$ only need to be $C^1$.

Theorem \ref{l1} also appears in {\cite{MR2597943,MR1046451}}. In
{\cite{MR1046451}}, Kulpa defines
\[ I (h^1, \ldots, h^n) = \int_{\Omega} \det h' (x) \mathd x \]
for smooth $h : \bar{\Omega} \rightarrow \mathbb{R}^n$ and manages to show
\begin{equation}
  I (f^1_+, f^2_+, \ldots, f_+^n) = I (f^1_-, f^2_+, \ldots, f_+^n) \label{e3}
\end{equation}
using Fubini theorem and projections $\Pi_i : \mathbb{R}^n \rightarrow
\mathbb{R}^{n - 1}$ given by deleting the $i$-th component. Then (\ref{e1})
follows by applying (\ref{e3}) successively to replace the remaining columns
of $\det f_+' (x)$ by those of $\det f_-' (x)$. While in {\citet[pp.\
464]{MR2597943}}, (\ref{e1}) is proved using the facts that $L (P) = \det P$
is a null Lagrangian and the energy $\int_{\Omega} L (w' (x), w (x), x) \mathd
x$ of null Lagrangian $L : \mathbb{M}^{m \times n} \times \mathbb{R}^m \times
\bar{\Omega} \rightarrow \mathbb{R}$ depends only on boundary conditions.

The purpose of this note is to present another two proofs of Theorem \ref{l1}.
The first one is performed by playing classical analysis, which is written for
readers not familiar with differential forms and in our opinion is very
elementary and transparent. The second one is quite short, which depends on
differential form and Stokes formula on manifolds. As demonstrated in the
paragraph after Theorem \ref{l1}, our note is useful for understanding the
Brouwer fixed point theorem from the analytic point of view. In addition to
{\cite{MR2597943,MR4725395,MR1046451,MR3711061}} mentioned above, other
analytic proofs of the Brouwer fixed point theorem can be found in \cite[pp.\ 467--470]{MR117523} and
{\cite{MR610487,MR1699248,MR505523,MR600910}}.

\section{Proof of Theorem \ref{l1}}\label{s2}

Let $f^i_{\pm}$ be the $i$-th component of $f_{\pm}$ and $\Phi_{\pm} =
(f_{\pm}^2, \ldots, f^n_{\pm})$. Then $\Phi_+ = \Phi_-$ on $\partial \Omega$.
Let
\[ A_{\pm} = \left( \frac{\partial (f_{\pm}^2, \ldots, f_{\pm}^n)}{\partial
   (x^2, \ldots, x^n)}, - \frac{\partial (f_{\pm}^2, f_{\pm}^3, \ldots,
   f_{\pm}^n)}{\partial (x^1, x^3, \ldots, x^n)}, \ldots, (- 1)^{n + 1}
   \frac{\partial (f_{\pm}^2, \ldots, f_{\pm}^n)}{\partial (x^1, \ldots, x^{n
   - 1})} \right) \text{,} \]
then the $i $-th component of $A_{\pm}$ is just the cofactor of
$\partial_{x^i} f_{\pm}^1$ in $\det f_{\pm}' (x)$. By Jacobi identity (see
e.g.\ {\citet[Eq.\ (5)]{MR610487}} or {\citet[Eq.\ (0.3)]{MR4701455}}), we have
$\mathrm{div} A_{\pm} = 0$. Thus
\begin{equation}
  \mathrm{div} (f_{\pm}^1 A_{\pm}) = \nabla f_{\pm}^1 \cdummy A_{\pm} +
  f_{\pm}^1 \mathrm{div} A_{\pm} = \nabla f_{\pm}^1 \cdummy A_{\pm} = \det
  f_{\pm}' (x) \text{.} \label{d}
\end{equation}
For $a \in \partial \Omega$, let $\eta : (u^1, \ldots, u^{n - 1}) \mapsto x$
be a local parametrization of $\partial \Omega$ at $a = \eta (\alpha)$ such
that
\[ N(a) = \left( \frac{\partial (x^2, \ldots, x^n)}{\partial (u^1, \ldots, u^{n -
   1})}, - \frac{\partial (x^1, x^3, \ldots, x^n)}{\partial (u^1, \ldots, u^{n
   - 1})}, \ldots, (- 1)^{n + 1} \frac{\partial (x^1, \ldots, x^{n -
   1})}{\partial (u^1, \ldots, u^{n - 1})} \right)_{\alpha} \]
is an outward normal vector of $\partial \Omega$ at $a$. Applying the chain
rule to the composition $y_{\pm} = \Phi_{\pm} (\eta (u))$, we see that
\[ \left( \begin{array}{ccc}
     \dfrac{\partial y_{\pm}^2}{\partial u^1} & \cdots & \dfrac{\partial
     y_{\pm}^2}{\partial u^{n - 1}}\\
     &  & \\
     \dfrac{\partial y_{\pm}^n}{\partial u^1} & \cdots & \dfrac{\partial
     y_{\pm}^n}{\partial u^{n - 1}}
   \end{array} \right)_{\alpha} = \left( \begin{array}{ccc}
     \dfrac{\partial y_{\pm}^2}{\partial x^1} & \cdots & \dfrac{\partial
     y_{\pm}^2}{\partial x^n}\\
     &  & \\
     \dfrac{\partial y_{\pm}^n}{\partial x^1} & \cdots & \dfrac{\partial
     y_{\pm}^n}{\partial x^n}
   \end{array} \right)_a \left( \begin{array}{ccc}
     \dfrac{\partial x^1}{\partial u^1} & \cdots & \dfrac{\partial
     x^1}{\partial u^{n - 1}}\\
     &  & \\
     \dfrac{\partial x^n}{\partial u^1} & \cdots & \dfrac{\partial
     x^n}{\partial u^{n - 1}}
   \end{array} \right)_{\alpha} \text{,} \]
where the matrix on the left is the Jacobian matrix $(\Phi_{\pm} \circ \eta)'
(\alpha)$.

Because $\Phi_+ = \Phi_-$ on $\partial \Omega$ and $\eta (u) \in \partial
\Omega$ for $u$ near $\alpha$, we have $\Phi_+ \circ \eta = \Phi_- \circ \eta$
near $\alpha$. Using the Cauchy-Binet formula we deduce
\begin{align*}
  (A_+ \cdummy N) (a) & =  \sum_{i = 1}^n \left. \frac{\partial (f_+^2,
  \ldots, f_+^n)}{\partial (x^1, \ldots, \hat{x}^i, \ldots, x^n)} \right|_a
  \cdummy \left. \frac{\partial (x^1, \ldots, \hat{x}^i, \ldots,
  x^n)}{\partial (u^1, \ldots, u^{n - 1})} \right|_{\alpha}\\
  & =  \det (\Phi_+ \circ \eta)' (\alpha) = \det (\Phi_- \circ \eta)'
  (\alpha) = (A_- \cdummy N) (a) \text{.}
\end{align*}
Therefore
\begin{equation}
  A_+ \cdummy \nu = \frac{A_+ \cdummy N}{| N |} = \frac{A_- \cdummy N}{| N |}
  = A_- \cdummy \nu \label{e5}
\end{equation}
on $\partial \Omega$, where $\nu = N / | N |$ is outward unit normal vector
field on $\partial \Omega$.

Since $f^1_+ = f^1_-$ on $\partial \Omega$, using the divergence theorem we
deduce from (\ref{d}) and (\ref{e5})
\begin{align*}
  \int_{\Omega} \det f_+' (x) \mathd x & =  \int_{\Omega} \mathrm{div} (f_+^1
  A_+) \mathd x = \int_{\partial \Omega} f_+^1 A_+ \cdummy \nu \mathd \sigma\\
  & =  \int_{\partial \Omega} f_-^1 A_- \cdummy \nu \mathd \sigma =
  \int_{\Omega} \det f_-' (x) \mathd x \text{.}
\end{align*}
\begin{example}
  Let $A$ be an $n \times n$ matrix, $\varphi \in C_0^{\infty} (\Omega,
  \mathbb{R}^n)$. Set $f_+ (x) = A x + \varphi (x)$, $f_- (x) = A x$, we get
  the formula stated in {\citet[pp. 21]{MR3887613}}:
  \[ \int_{\Omega} \det (A + \varphi' (x)) \mathd x = \det (A) \cdummy m
     (\Omega) \text{.} \]
\end{example}

\section{Another proof using differential forms}

To conclude this note, we present another proof using differential forms. Let
$f^i$ be the components of $f:\bar{\Omega}\to\mathbb{R}^n$, then it is well known that
\begin{align*}
  \det f' (x) \mathd x^1 \wedge \cdots \wedge \mathd x^n & =  \mathd f^1
  \wedge \mathd f^2 \wedge \cdots \wedge \mathd f^n \text{.}\\
  & =  \mathd (f^1 \mathd f^2 \wedge \cdots \wedge \mathd f^n) \text{,}
\end{align*}
where in the second equality we used $\mathd^2 = 0$. Recall the Stokes formula
\[ \int_{\Omega} \mathd \omega = \int_{\partial \Omega} i^{\ast} \omega
   \text{,} \]
where $\omega$ is an $(n - 1)$-form on $\Omega$ and $i^{\ast}$ is the pull
back induced by the embedding $i : \partial \Omega \rightarrow \Omega$.
Because $i^{\ast} \circ \mathd = \mathd \circ i^{\ast}$, we conclude
\begin{align*}
  \int_{\Omega} \det f' (x) \mathd x & =  \int_{\Omega} \det f' (x) \mathd
  x^1 \wedge \cdots \wedge \mathd x^n\\
  & =  \int_{\Omega} \mathd (f^1 \mathd f^2 \wedge \cdots \wedge \mathd f^n)
  = \int_{\partial \Omega} i^{\ast} (f^1 \mathd f^2 \wedge \cdots \wedge
  \mathd f^n)\\
  & =  \int_{\partial \Omega} (f^1 \circ i) \mathd (f^2 \circ i) \wedge
  \cdots \wedge \mathd (f^n \circ i) \text{.}
\end{align*}
Clearly the right hand side depends only on $f \circ i = f |_{\partial \Omega}
 $.

\end{document}